\documentclass[12pt]{amsart}

\usepackage{amssymb}
\usepackage{pb-diagram}
\usepackage{enumitem} 
\usepackage[dvipsnames]{xcolor}
\usepackage{tikz-cd}
\usepackage{graphicx}
\usepackage{amsmath}
\numberwithin{equation}{section}
\usepackage[all]{xy}

\newtheorem{Theorem}{Theorem}[section]
\newtheorem{Proposition}[Theorem]{Proposition}
\newtheorem{Lemma}[Theorem]{Lemma}
\newtheorem{Corollary}[Theorem]{Corollary}

\newtheorem{Remark}[Theorem]{Remark}

\numberwithin{equation}{section}

\newcommand{\Hom}{\text{Hom}}
\newcommand{\Ext}{\text{Ext}}
\newcommand{\Quot}{\text{Quot}}
\newcommand{\Cliff}{\text{Cliff}}
\newcommand{\OO}{\mathcal O}
\newcommand{\OS}{\mathcal O_S}
\newcommand{\CC}{\mathbb{C}}
\newcommand{\EE}{\mathcal{E}}
\newcommand{\MM}{\mathcal{M}}
\newcommand{\red}{\mathrm{red}}
\newcommand{\PP}{\mathcal{P}}
\newcommand{\GG}{\mathcal{G}}
\newcommand{\WW}{\mathcal{W}}
\DeclareMathOperator{\Coh}{Coh}
\DeclareMathOperator{\Pic}{Pic}
\DeclareMathOperator{\NS}{NS}
\DeclareMathOperator{\ch}{ch}
\DeclareMathOperator{\rk}{rk}
\newcommand{\WallCat}{\mathcal P_0}

\begin{document}

\baselineskip=16pt

\title{Stability conditions on Lazarsfeld-Mukai bundles over abelian surfaces}

\author{Erick Luna}
\address{Universidad Autónoma de Aguascalientes
\newline  Av. Universidad \# 940, Ciudad Universitaria, \newline C.P. 20100, Aguascalientes, Aguascalientes, M\'exico.}
\email{erick.luna@edu.uaa.mx}

\author{Graciela Reyes-Ahumada}
\address{SECIHTI-ITCM
\newline 1° Mayo \& Sor Juana Inés de la Cruz, Los Mangos.  \newline C.P. 89440 Cd. Madero, Tam. México}
\email{grace@cimat.mx}



\date{\today}

\begin{abstract}
We study the stability of Lazarsfeld-Mukai bundles associated with smooth curves over abelian surfaces. When the general curve in $|L|$ have Clifford dimension one and maximal gonality, for $\rho \gg 0$ we show that every component of the relative Brill-Noether loci \(W^1_d(|L|)\) dominating the linear system parametrizes generically \(\mu_L\)-stable rank-two Lazarsfeld-Mukai bundles. We also apply Bayer-Li wall-crossing to higher-rank bundles: if \(\operatorname{NS}(S)=\mathbb Z\cdot L\), \(A\) is a globally generated complete linear series of degree \(d\) on a smooth curve \(C\in|L|\), \(h^0(C,A)=r+1\) with \(r>2\), and \(d\ne g-1\), then the associated bundle \(E_{C,A}\) is \(\mu_L\)-stable and \(L\)-Gieseker-stable. The corresponding object \(E_{C,A}^{\vee}[1]\) for \(d<g-1\), or \(E_{C,A}\) for \(d>g-1\), is stable in the adjacent chamber, as its reduction. These results connect dimension estimates for Brill-Noether loci with stability and wall-crossing on abelian surfaces.
\end{abstract}

\maketitle

\section{Introduction}

Let $S$ be a smooth projective abelian surface, let $L$ be ample, and write $|L|_s$ for its smooth connected curves, each of genus $g=1+L^2/2$. For such a curve $C$, set $W_d^1(C)=\{A\in\Pic^d(C):h^0(C,A)\geq2\}$; the relative locus $\mathcal W_d^1(|L|_s)$ parametrizes pairs $(C,A)$ with $A\in W_d^1(C)$. Its expected fiber dimension, when nonnegative, is the Brill-Noether number $\rho(g,1,d)=g-2(g-d+1)=2d-g-2$, so the expected relative dimension is $g-2+\rho(g,1,d)$; a negative Brill-Noether number predicts an empty general fiber. The gonality is the least degree of a pencil, and the rank-two discussion assumes general curves have Clifford dimension one (their Clifford index is computed by a pencil) and maximal gonality $k=\lfloor(g+3)/2\rfloor$. Knutsen, Lelli-Chiesa, and Mongardi study Severi varieties and Brill-Noether loci on general polarized abelian surfaces \cite{knutsen2019severi}, while Bayer-Li describe Brill-Noether theory for curves in primitive systems on generic abelian surfaces \cite{bayer2018brill}. Moretti's recent work computes gonality and the dimensions of $W_d^1(C)$ for general curves in non-primitive systems on simple abelian surfaces, and derives Green's conjecture in that setting \cite{moretti2025brill}.

The Lazarsfeld-Mukai construction associates a vector bundle on $S$ to a series on $C$ by extending its evaluation data from the curve to the surface; a complete base-point-free $\mathfrak{g}_d^r$ gives a bundle of rank $r+1$. Its stability connects the geometry of this series with moduli spaces of sheaves: the slope stability compares $\mu_L(E)=c_1(E)\cdot L/\operatorname{rk}(E)$ with subsheaf slopes, while Gieseker stability is defined using reduced Hilbert polynomials. Lazarsfeld's foundational work on Brill-Noether-Petri theory for curves on $K3$ surfaces established a vector-bundle approach that informs later Lazarsfeld-Mukai methods \cite{lazarsfeld1986brill}. Narayanan studies semistability of Lazarsfeld-Mukai bundles on abelian surfaces and stability for dominating components on associated Kummer surfaces \cite{narayanan2018semistability}.

For a complete base-point-free pencil on a smooth curve in an ample linear system, the associated LM bundle has rank two. We study possible destabilizing filtrations and use intersection-theoretic and parameter-space estimates to assess how non-stable bundles may contribute to the relative Brill-Noether locus. The Clifford-dimension and gonality hypotheses enter these estimates; a bound on the non-stable sublocus alone does not establish the expected dimension or reducedness of the full Brill-Noether locus.

The final part applies Bayer-Li wall-crossing to complete, globally generated series of higher dimension \cite{bayer2018brill}. If $S$ has $\NS(S)=\mathbb Z\cdot L$, $C\in|L|$ is smooth, and $A$ is a globally generated line bundle of degree $d$ with $h^0(C,A)=r+1$, for $r>2$, and $d\ne g-1$, then $E_{C,A}$ is $\mu_L$-stable and $L$-Gieseker-stable. Wall-crossing varies the stability condition across a wall in the derived category, where stability can change; this higher-rank argument is independent of the rank-two estimates. Bayer explains how Bridgeland wall-crossing recovers Lazarsfeld's Brill-Noether theorem for curves on $K3$ surfaces and surveys applications of stability conditions on surfaces \cite{bayer2016wall}. Feyzbakhsh's effective restriction theorem gives a complementary criterion for stability of restrictions of sheaves to divisors \cite{feyzbakhsh2022effective}.

Section 2 recalls the LM construction and its basic invariants; Section 3 studies rank-two stability and parameter estimates; Section~4 introduces the Bridgeland wall, states the Bayer-Li results used here, and proves the higher-rank stability theorem.

\section{Lazarsfeld-Mukai bundles for abelian surfaces}

Let \(A\) be a complex abelian surface. For any coherent sheaf \(Q\) on \(A\), the Euler characteristic is defined as the alternating sum of the cohomology dimensions
\[
\chi(A,Q) \colon = \sum_{i=0}^{2} (-1)^i \dim H^{i}(A,Q).
\]
For \(A\) an abelian surface, the tangent bundle has vanishing Chern classes \(c_1(T_A) = 0\) and \(c_2(T_A) = 0\), so the Todd class \(\operatorname{td}(A)\) is trivial. Following \cite{huybrechts2005complex}, the Hirzebruch-Riemann-Roch theorem then yields
\[
\chi(A,Q) = \int_A \operatorname{ch}(Q),
\]
where \(\operatorname{ch}(Q)\) is the Chern character of \(Q\). In terms of the Chern classes \(r = \operatorname{rk}(Q)\), \(c_1 = c_1(Q)\), and \(c_2 = c_2(Q)\), this becomes
\[
\chi(A,Q) = \frac{1}{2}\, c_1^2 - c_2,
\]
with all integrals evaluated against the fundamental class of \(A\).
In the special case where \(Q = L\) is a line bundle, the rank \(r\) is 1 and \(c_2(L) = 0\), so the Euler characteristic simplifies to
\[
\chi(A,L) = \frac{1}{2}\, c_1(L)^2.
\]

\noindent As in \cite{huybrechts2010geometry}, we have that the Mukai vector of a coherent sheaf \(Q\) on a surface \(A\) is defined by applying the Chern character to the modified Todd class, namely
\[
v(Q) \colon = \operatorname{ch}(Q)\sqrt{\operatorname{td}(A)}.
\]
On an abelian surface \(A\), the Todd class is trivial, so \(v(Q) = \operatorname{ch}(Q)\) lies in
\[
H^0(A,\mathbb{Z}) \oplus H^2(A,\mathbb{Z}) \oplus H^4(A,\mathbb{Z}).
\]
In terms of the Chern classes \(r = \operatorname{rk}(Q)\), \(c_1 = c_1(Q)\), and \(c_2 = c_2(Q)\), the Mukai vector reads
\[
v(Q) = \big(r,\ c_1(Q),\ \tfrac{1}{2}\, c_1(Q)^2 - c_2(Q)\big).
\]
Derived equivalences and Fourier-Mukai transforms transport this vector by the induced cohomological isometry, and it is a fundamental invariant in the study of moduli spaces of sheaves on abelian surfaces.
If \(Q = L\) is a line bundle on \(A\), then \(r = 1\) and \(c_2(L) = 0\), so the Mukai vector becomes
\[
v(L) = \big(1,\ c_1(L),\ \tfrac{1}{2}\, c_1(L)^2\big).
\]

Let \(A\) be a complex abelian surface and let \(Q\) and \(R\) be coherent sheaves on \(A\), with Mukai vectors
\[
v(Q) = \big(v_0(Q),\, v_1(Q),\, v_2(Q)\big), \quad
v(R) = \big(v_0(R),\, v_1(R),\, v_2(R)\big),
\]
lying in \(H^0(A,\mathbb{Z}) \oplus H^{2}(A,\mathbb{Z}) \oplus H^{4}(A,\mathbb{Z})\). The Mukai pairing of $v(Q)$ and $v(R)$ is the symmetric bilinear form defined as
\[
\langle v(Q), v(R)\rangle \colon = -\sum_{i=0}^{2} (-1)^i \dim \operatorname{Ext}^{i}_A(Q,R).
\]
If one writes \(v(Q) = (r, c_1, v_2)\), so \(v_0 = r\), \(v_1 = c_1\), \(v_2 = \frac{1}{2}c_1^2 - c_2\), and similarly \(v(R) = (r', c_1', v_2')\), then the Mukai pairing reads explicitly
\[
\langle v(Q), v(R)\rangle
= v_1(Q)\cdot v_1(R) - v_0(Q)v_2(R) - v_2(Q)v_0(R),
\]
where \(v_1(Q)\cdot v_1(R)\) denotes the intersection product \(\int_A v_1(Q)\wedge v_1(R)\).

Let \(A\) be a complex abelian surface and let \(M\) be a coherent sheaf on \(A\) with fixed Mukai vector \(v(M) = (r, c_1, v_2)\). Fix an ample divisor class \(H\) on \(A\) and consider the moduli space \(\mathcal{M}_H^\mathrm{ss}(v(M))\) of \(H\)-semistable sheaves on \(A\) with Mukai vector \(v(M)\), in the sense of geometric invariant theory as developed in \cite{huybrechts2010geometry}. The open subscheme \(\mathcal{M}_H^\mathrm{s}(v(M)) \subset \mathcal{M}_H^\mathrm{ss}(v(M))\) parametrizes \(H\)-stable sheaves and, when this moduli space is nonempty, \(v(M)\) is primitive, and \(H\) is generic for \(v(M)\), stable and semistable objects coincide and the moduli space is smooth, irreducible, and projective over \(\mathbb{C}\). 

\noindent The expected dimension of the moduli space is governed by the deformation theory of the sheaf. For a stable sheaf \(E\) with \(\operatorname{rk}(E)=r\) and Mukai vector \(v(E)\), the tangent space at \([E]\) is identified with \(\operatorname{Ext}^1(E,E)\) and the obstruction space with \(\operatorname{Ext}^2(E,E)\). On an abelian surface, the trace-free obstruction space vanishes, whereas $\Ext^2(E,E)\simeq\mathbb C$; the stable coarse moduli space therefore has dimension \(\langle v(E),v(E)\rangle+2\); see \cite{huybrechts2010geometry}. 

When \(M\) is a line bundle, \(v(M) = (1, c_1(M), \tfrac{1}{2}c_1(M)^2)\) and the Mukai self-pairing vanishes, hence the stable moduli space has dimension \(2\).

Let $S$ be a smooth complex projective abelian surface, let
$i:C\hookrightarrow S$ be a smooth connected curve, and denote as
$L=\mathcal O_S(C)$ and $g=1+L^2/2$. Let $A$ be a globally generated
line bundle of degree $d$ on $C$, with $h^0(C,A)=n=r+1$.
Thus we use a complete $g_d^r$, not an arbitrary point of $W_d^r(C)$
with possibly more sections. Evaluation on $S$ is surjective
\begin{equation}\label{eq:1}
 0\to F_{C,A}\to H^0(C,A)\otimes\OS
 \xrightarrow{\operatorname{ev}_{A,S}}i_*A\to0,
\end{equation}
and the kernel is locally free. The Lazarsfeld-Mukai bundle associated to
$(C,A)$ is $E_{C,A}=F_{C,A}^\vee$. Adjunction gives
$\omega_C\simeq L|_C$, and dualizing \eqref{eq:1} yields
\begin{equation}\label{eq:2}
 0\to H^0(C,A)^\vee\otimes\OS\to E_{C,A}
 \to i_*(\omega_C\otimes A^{-1})\to0.
\end{equation}
Both quotients are pushforwards from $C$, not vector bundles on $S$.
In particular, $E_{C,A}$ has a distinguished $n$-dimensional subspace
of sections; this need not be its entire space of sections.

\begin{Theorem}[Basic LM invariants]\label{thm:basic-lm-invariants}
With the above assumptions, write $E=E_{C,A}$ and
$B=\omega_C\otimes A^{-1}$. Then:
\begin{itemize}
\item $\rk(E)=n$, $\det(E)=L$, $c_2(E)=d$, and
$\chi(S,E)=g-1-d$;
\item $h^2(S,E)=0$. For the map
$\gamma:H^0(C,A)\otimes H^1(S,\OS)\to H^1(C,A)$ induced by evaluation,
\[
 \tau:=h^1(S,E)=\dim\ker\gamma,\qquad 0\leq\tau\leq2n;
\]
\item $h^0(S,E)=g-1-d+\tau
=n+h^1(C,A)-\operatorname{rank}\gamma\geq n$;
\item the distinguished sections generate $E$ away from $C$.
If the image of $H^0(S,E)\to H^0(C,B)$ generates $B$ away from a
finite subset of $C$, then $E$ is generated outside a finite subset of $S$;
\item $v(E)=(n,L,g-1-d)$ and
$v(E)^2=L^2-2n(g-1-d)$.
\end{itemize}
\end{Theorem}
\begin{proof}
The character of $i_*A$ is $(0,L,d+1-g)$. Additivity in
\eqref{eq:1} and duality give the rank, determinant, Chern classes
and Euler characteristic. The induced map
$H^0(C,A)\otimes H^0(S,\OS)\to H^0(C,A)$ is the identity.
Hence $H^0(S,F_{C,A})=0$ and
$H^1(S,F_{C,A})\simeq\ker\gamma$. Serre duality gives
$H^2(S,E)=0$ and $h^1(S,E)=h^1(S,F_{C,A})$.
Now $h^1(S,\OS)=2$ and Riemann-Roch prove the cohomology formulas.
The generation statements follow from \eqref{eq:2}; the Mukai square
uses the convention fixed above.
\end{proof}
\noindent Notice that from the non-vanishing cohomology $H^1(S,\OS)$  the K3 formulas cannot be transferred unchanged to an abelian surface. In particular, neither
$h^0(S,E)=n$ nor automatic generation along $C$ is asserted.
For the restriction and Petri calculations, define the kernel on $C$ by
\[
 0\to M_A\to H^0(C,A)\otimes\mathcal O_C\to A\to0.
\]
Restricting \eqref{eq:1} to $C$ and accounting for its Tor term gives
\begin{equation}\label{eq:3}
 0\to A\otimes\omega_C^{-1}\to F_{C,A}|_C\to M_A\to0.
\end{equation}
Equivalently, tensoring by $B$ gives
\[
 0\to\mathcal O_C\to F_{C,A}|_C\otimes B\to M_A\otimes B\to0.
\]
The Petri map
$\mu_A:H^0(C,A)\otimes H^0(C,B)\to H^0(C,\omega_C)$ has kernel
$H^0(C,M_A\otimes B)$. These are statements on the curve; they do not
supply an isomorphism involving $H^0(S,E_{C,A}\otimes F_{C,A}^\vee)$.

\begin{Proposition}[\cite{lelli2013stability}, Prop.~3.1]
Let $E$ and $F$ be torsion-free sheaves on $S$. Then:
\begin{enumerate}
\item If $E$ and $F$ are $\mu_H$-semistable with $\mu_H(E) > \mu_H(F)$, then $\Hom(E,F) = 0$.
\item If $E$ and $F$ are $\mu_H$-stable with $\mu_H(E) = \mu_H(F)$ and $0 \neq \phi \in \Hom(E,F)$, then $\operatorname{rk}(E) = \operatorname{rk}(F)$ and $\phi$ is injective and an isomorphism away from a finite subset of $S$.
\end{enumerate}
\end{Proposition}

\begin{Remark}\label{h2}
    On an abelian surface $S$ since $K_S \cong \mathcal{O}_S$ by Serre duality $H^2(S,E)\cong \Hom(E,\OS)^\vee$ and the first part of the above Proposition implies that, if $E$ is $\mu_H$-semistable and $\mu_H(E)>0$ then $h^2(S, E)=0$. Also if $E$ is $\mu_H$-semistable and $\mu_H(E)<0$ then $h^0(S, E)=0$. 
\end{Remark}

\begin{Lemma}[\cite{lelli2013stability}, Lemma~3.2] \label{Lemma 3.2}
Let $E, Q \in \Coh(S)$ be torsion-free sheaves with $\rk(E)\geq2$ and $\rk(Q)>0$. Suppose $E$ is globally generated outside a finite set of points, $h^2(S, E) = 0$ and there exists a surjective morphism $\varphi \colon E \to Q$. Then $h^0(S, Q^{\vee\vee}) \geq 2$.
\end{Lemma}

\begin{proof}
As a quotient of $E$, the sheaf $Q$ is globally generated outside a finite set. If $\rk(Q) \geq 2$, then
\[
h^0(S, Q^{\vee\vee}) \geq h^0(S, Q) \geq 2,
\]
as required.

\noindent Now suppose $\rk(Q)=1$. Write $Q=N\otimes I_\xi$ where
$N=Q^{\vee\vee}$ is a line bundle and $\xi$ is zero-dimensional.
Generation outside a finite subset implies that $N$ has a nonzero section
and no fixed divisorial component. If $h^0(S,N)=1$, its unique section
has no divisorial zero locus and therefore trivializes $N$.
The nonzero composite $E\twoheadrightarrow Q\hookrightarrow N\simeq\OS$
would contradict $h^2(S,E)=\dim\Hom(E,\OS)=0$.
Hence $h^0(S,N)\geq2$. No assertion of base-point-freeness is needed;
a line bundle with no fixed divisorial component may have isolated
base points.
\end{proof}

\section{Stability of Lazarsfeld-Mukai bundles of rank 2}
Let $S$ be a smooth, projective abelian surface and consider a line bundle $L\in Ample(S)$ such that a general curve $C\in |L|_s$ has genus $g$, Clifford dimension $1$ and maximal gonality $k =
\lfloor\frac{g+3}{2}\rfloor$. In this section we prove that, if $C$ is general in its linear system and $\rho(g, 1, d) > 0$, the LM bundle associated with a general complete, base-point-free $g^1_d$ on $C$ is $\mu_L$-stable.

\noindent Fix a rank-2 LM bundle $E = E_{C,A} $ corresponding to a complete, base-point-free pencil $A\in  W^1_d (C)$ with $C \in|L|_s$; we assume $E$ is not $\mu_L$-stable. In the case where $E$ is $\mu_L$-unstable (resp. properly $\mu_L$-semistable) we consider its HN filtration (resp. JH filtration) $0 \subset M \subset E$, which gives a short exact sequence
\begin{equation}\label{eqMN}
    0 \rightarrow M \rightarrow E \rightarrow N \otimes I_\xi \rightarrow 0,
\end{equation}
where $M$ and $N$ are two line bundles such that 
$$\mu_L(M) > \mu_L(E) = g-1 > \mu_L(N) \text{(resp.} \mu_L(M) = \mu_L(E) = \mu_L(N)\text{)}$$ 
and $I_\xi$ is the ideal sheaf of a 0-dimensional subscheme $\xi \subset S$ of length $l = d - c_1(N) \cdot c_1(M)$. 

\begin{Lemma}
If general curves in $|L|_s$ have Clifford dimension 1 and constant gonality $k$, one has $c_1(M ) \cdot c_1(N ) \geq k-4$
\end{Lemma} 

 \begin{proof}
Since $g-1<\mu_L(M)=c_1(M)^2+c_1(N)\cdot c_1(M)$ then if $c_1(M)^2\leq 0$ we get
$$g\leq c_1(M)^2+c_1(N)\cdot c_1(M)\leq c_1(N)\cdot c_1(M).$$
 On the other hand, from hypothesis in the gonality 
 $$k-1\leq \frac{g+1}{2}\leq g.$$
 By putting together $k-1\leq c_1(N)\cdot c_1(M)$ and the lemma holds. From now on we assume  $$c_1(M)^2>0.$$
 Since $\mu_L(M)>0$ from Remark~\ref{h2} it follows $h^2(S,M)=0$ and $\chi(M)=h^0(S,M)-h^1(S,M)$.
 
 \begin{description}
    \item[Case $h^0(S,M)=0$] From Riemann-Roch theorem $\chi(M)=-h^1(S,M)=c_1(M)^2/2$ which gives a contradiction on the assumption.
   
    \item[Case $h^0(S,M)=1$] In this case $\chi(M)=1-h^1(S,M)=c_1(M)^2/2$ from our assumption $h^1(S,M)=0$ and $c_1(M)^2=2$ then
    $$g-1<\mu_L(M)=c_1(M)^2+c_1(N)\cdot c_1(M)=2+c_1(N)\cdot c_1(M),$$
      and $g-3<c_1(N)\cdot c_1(M)$. From the hypothesis in the gonality
        $$k-4\leq g-3,$$
    it follows that $k-4< c_1(N)\cdot c_1(M)$.
    
    \item[Case $h^0(S,M)\geq 2$] 
   From Lemma \ref{Lemma 3.2}, we know that $h^0(S,N)\geq 2$ then $c_1(N)$ is effective and by Nakai-Moishezon criterion since $L$ is ample $0<c_1(L)\cdot c_1(N)= \mu_L(N)$ then Remark \ref{h2} implies $h^2(S,N)=0$.
     Since $M=L\otimes N^\vee$ then
$0=h^2(S,N)=h^0(S,N^\vee)=h^0(S,M(-C))$. By adjunction formula over $C$ we get
    $$M|_C=(L\otimes N^\vee)|_C\cong \omega_C \otimes N^\vee|_C,$$ and from the  cohomology of 
$$0\rightarrow M(-C) \rightarrow M\rightarrow M|_C\rightarrow 0,$$
it follows that $2\leq h^0(S,M)\leq h^0(C,M|_C)=h^0(C,\omega_C \otimes N^\vee|_
    C)=h^1(C,N|_C)$.
    Similarly the exact sequence
$$0\rightarrow M^\vee \rightarrow N\rightarrow N|_C\rightarrow 0,$$
gives $2\leq h^0(S,N)\leq h^0(C,N|_C),$
then $N|_C$ contributes to $\Cliff(C)$,
\begin{align*}
    \Cliff(N|_C) &= c_1(N)\cdot(c_1(N)+c_1(M))-2h^0(N|_C)+2\\
    &=c_1(N)^2+c_1(N)\cdot c_1(M)-2h^0(N|_C)+2
\end{align*}
since $h^2(S,N)=0$ then
$h^0(C,N|_C)\geq h^0(S,N)= \chi(N)+h^1(S,N)$ and
\begin{align*}
    \Cliff(N|_C) &\leq c_1(N)^2+c_1(N)\cdot c_1(M)-2\chi(N)-2h^1(S,N)+2\\
   &=c_1(N)^2+c_1(N)\cdot c_1(M)-2(c_1(N)^2/2)-2h^1(S,N)+2 \\
   &= c_1(N)\cdot c_1(M)-2h^1(S,N)+2
\end{align*}

From $\Cliff(N|_C)\geq k-2$ follows that
$$c_1(N)\cdot c_1(M)\geq k-4+2h^1(S,N)\geq k-4.$$
\end{description}
\end{proof}

As in \cite{lelli2013stability}, we use the following notation:
Fix a nonnegative integer $l$ and a non-trivial, globally generated line bundle $N$ on $S$ such that, having defined $M= L \otimes N^\vee$ either $\mu_L(M)=\mu_L(N) = g-1$ or
$\mu_L(M ) > g-1 > \mu_L(N)$. We consider the moduli stack $\EE_{N,l}$ parametrizing filtrations
$0 \subset M \subset E$ with $[M] \in \MM(v(M))(\CC)$ and $[E/M] \in \MM(v(N \otimes I_\xi))(\CC)$, where $l(\xi) = l$.

\begin{Remark}\label{morfp}
The dimension $\dim\EE_{N,l}$ can be computed by considering the natural morphism 
$$p: \EE_{N,l} \rightarrow \MM(v(M))\times \MM(v(N\otimes I_\xi))$$ sending a short exact sequence as in \eqref{eqMN} to $(M,N\otimes I_\xi)$. The fiber of $p$ at a $\CC$-point $(M,N\otimes I_\xi)$ is the quotient stack $[\Ext^1(N\otimes I_\xi, M)/ \Hom(N\otimes I_\xi, M)]$ with trivial action. The dimension of the fiber is $\dim\Ext^1(N\otimes I_\xi, M)-\dim\Hom(N\otimes I_\xi, M)$.
From the duality on abelian varieties $\Ext^i(F, G)\cong \Ext^{2-i}(G, F)^\vee$  then
\begin{align*}
 \dim\Ext^1(N\otimes I_\xi, M)  &= \dim\Ext^0(N\otimes I_\xi, M)+\dim\Ext^2(N\otimes I_\xi, M)-\chi(N\otimes I_\xi,M)\\
&= \dim\Ext^0(N\otimes I_\xi, M)+\dim\Ext^0(M,N\otimes I_\xi)-\chi(N\otimes I_\xi,M)
\end{align*}
since $\Ext^0(F, G)\cong\Hom(F,G)$ and $\chi(F,G)=-\langle v(F), v(G)\rangle$ the dimension of the fiber is
$$\dim\Hom(M,N\otimes I_\xi)-\chi(N\otimes I_\xi,M)=\dim\Hom(M,N\otimes I_\xi)+\langle v(N\otimes I_\xi),v(M)\rangle,$$
we conclude that
\[
 \dim\EE_{N,l}=2l+\dim\Hom(M,N\otimes I_\xi)
 +\langle v(N\otimes I_\xi),v(M)\rangle.
\]\end{Remark} 

\begin{Remark}
Let \(M\) and \(N\) be line bundles on the surface, and let \(\xi\) be a finite subscheme of length \(l\), with associated ideal sheaf \(I_\xi\). We denote by \(N\otimes I_\xi\) the sheaf obtained by tensoring \(N\) with \(I_\xi\). The goal is to compute the Mukai pairings, using the componentwise formula
\[
\langle v,w\rangle=v_1\cdot w_1-v_0w_2-v_2w_0.
\]
The Mukai vector of a line bundle \(M\) is
\[
 v(M)=\big(1,c_1(M),\tfrac12c_1(M)^2\big).
\]
Similarly, for \(N\) one has
\[
 v(N)=\big(1,c_1(N),\tfrac12c_1(N)^2\big).
\]
To compute the Mukai vector of \(N\otimes I_\xi\), consider the short exact sequence
\[
 0\to N\otimes I_\xi\to N\to N|_\xi\to0.
\]
The sheaf $N|_\xi$ has character $(0,0,l)$. Thus
\begingroup
\begin{align*}
 \ch(N\otimes I_\xi)
 &=\ch(N)-\ch(N|_\xi)\\
 &=\big(1,c_1(N),\tfrac12c_1(N)^2\big)-(0,0,l)\\
 &=\big(1,c_1(N),\tfrac12c_1(N)^2-l\big).
\end{align*}
\endgroup
Putting \(v=v(M)\), the self-pairing is
\begingroup
\begin{align*}
 \langle v(M),v(M)\rangle
 &=c_1(M)^2-\frac{c_1(M)^2}{2}-\frac{c_1(M)^2}{2}=0.
\end{align*}
\endgroup
Thus the Mukai self-pairing of a line bundle is
\[
 \langle v(M),v(M)\rangle=0.
\]
Let \(v=v(N\otimes I_\xi)\). Its rank remains one, not $1-l$. Then
\begingroup
\begin{align*}
 \langle v(N\otimes I_\xi),v(N\otimes I_\xi)\rangle
 &=c_1(N)^2-2\left(\frac{c_1(N)^2}{2}-l\right)=2l.
\end{align*}
\endgroup
Therefore the self-pairing of the twisted sheaf is
\[
 \langle v(N\otimes I_\xi),v(N\otimes I_\xi)\rangle=2l.
\]
Now take
\begingroup
\[
 v(M)=(1,c_1,\tfrac12c_1^2),\qquad
 v(N\otimes I_\xi)=(1,c_1',\tfrac12(c_1')^2-l),
\]
\endgroup
with \(c_1=c_1(M)\), \(c_1'=c_1(N)\). The Mukai pairing is
\begingroup
\begin{align*}
 \langle v(M),v(N\otimes I_\xi)\rangle
 &=c_1(M)\cdot c_1(N)
   -\left(\frac{c_1(N)^2}{2}-l\right)-\frac{c_1(M)^2}{2}\\
 &=c_1(M)\cdot c_1(N)-\frac{c_1(M)^2+c_1(N)^2}{2}+l.
\end{align*}
\endgroup
In terms of Chern classes, this is
\begin{equation}\label{eq:initial-rank-one-vectors}
\langle v(M),v(N\otimes I_\xi)\rangle
 =c_1(M)\cdot c_1(N)-\frac{c_1(M)^2+c_1(N)^2}{2}+l.
\end{equation}
In particular, when \(l=0\), this recovers \(\langle v(M),v(N)\rangle\).
\begingroup
For these stable rank-one sheaves, the coarse moduli dimensions are
\[
 \dim M_H(v(M))=2,\qquad
 \dim M_H(v(N\otimes I_\xi))=2l+2.
\]
Their moduli stacks have dimensions $1$ and $2l+1$, respectively.
These dimensions follow from $v^2+2$ for coarse spaces and $v^2+1$
for stable stacks.
\endgroup
\end{Remark}

\begin{Remark}\label{morfq}
Consider the natural morphism  
$$q: \EE_{N,l} \rightarrow \MM(v(E))$$ mapping an exact sequence as \ref{eqMN} to $[E]$. Denote as $\tilde{\PP}_{N,l}=\overline{\operatorname{Im}(q)}$ the closure of the image of $q$ and consider the open substack $\PP_{N,l}$ inside  $\tilde{\PP}_{N,l}$ whose points correspond to vector bundles $E$ satisfying $h^2(S, E)=0$ and $h^1(S,E)=\tau$, with $\tau$ fixed. This last condition defines a locally closed cohomology stratum, not necessarily an open substack; we use $\PP_{N,l}$ for that stratum. 

\noindent Consider the Grassmann bundle $\GG_{N,l}\rightarrow \PP_{N,l}$ with fiber over a point $[E]$ equal to $G(2, H^0(S, E))$. By Riemann-Roch and $h^2(S,E)=0$, $h^0(S,E)=g-1-d+\tau$, so $\dim G(2,H^0(S,E))=2(g-3-d+\tau)$ where $g-1-d+\tau\geq2$ and if $\PP_{N,l}$ is not empty then
$$\dim\GG_{N,l}=\dim\PP_{N,l}+2(g-3-d+\tau).$$
\end{Remark}

In order to estimate $\tau$ by the cohomology of the filtration. We use the long exact cohomology sequences of \eqref{eqMN} and
$0\to N\otimes I_\xi\to N\to N|_\xi\to0$.
Put $Q=N\otimes I_\xi$. Since $\mu_L(M)\geq g-1>0$,
Serre duality gives $H^2(S,M)=0$. For the connecting map
$\delta:H^0(S,Q)\to H^1(S,M)$, the first sequence gives
\[
 0\to\operatorname{coker}\delta\to H^1(S,E)\to H^1(S,Q)\to0.
\]
Define the rank of restriction to the length-$l$ subscheme by
\[
 r_\xi=\operatorname{rank}\bigl(H^0(S,N)\to H^0(\xi,N|_\xi)\bigr).
\]
The second sequence, together with $h^0(\xi,N|_\xi)=l$ and the
vanishing of higher cohomology of $N|_\xi$, gives
\[
 h^0(S,Q)=h^0(S,N)-r_\xi,\qquad
 h^1(S,Q)=h^1(S,N)+l-r_\xi.
\]
Consequently the chosen extension satisfies the exact identity
\begin{equation}\label{eq:tau-filtration-identity}
 \tau=h^1(S,M)+h^1(S,N)+l-r_\xi-\operatorname{rank}\delta.
\end{equation}
Using $0\leq\operatorname{rank}\delta\leq h^1(S,M)$ yields
\[
 h^1(S,N)+l-r_\xi\leq\tau
 \leq h^1(S,M)+h^1(S,N)+l-r_\xi.
\]
For a complete LM pencil, $\tau\leq4$ follows from the original
evaluation sequence: $H^1(S,E^\vee)$ embeds into
$H^0(C,A)\otimes H^1(S,\OS)$, a four-dimensional space, and
Serre duality identifies its dimension with $h^1(S,E)$.
Also $h^2(S,E)=0$, $\chi(S,E)=g-1-d$, and $h^0(S,E)\geq2$
imply $\tau\geq d+3-g$. Thus the estimate used below is
\begin{equation}\label{eq:tau-filtration-bound}
 \tau\leq\min\{4,h^1(S,M)+h^1(S,N)+l-r_\xi\},
\end{equation}
with lower bound
$\tau\geq\max\{0,d+3-g,h^1(S,N)+l-r_\xi\}$.

\noindent This bound uses only the filtration and completeness of the pencil;
it does not assume that both endpoints are ample. Since $N$ is
nontrivial and globally generated, $h^2(S,N)=0$. Hence the endpoint
terms can be computed by Riemann-Roch:
\[
 h^1(S,M)=h^0(S,M)-\tfrac12c_1(M)^2,\qquad
 h^1(S,N)=h^0(S,N)-\tfrac12c_1(N)^2.
\]
The remaining defect $l-r_\xi$ records the failure of $\xi$ to impose
independent conditions on sections of $N$. Global generation gives
$r_\xi\geq1$ if $l>0$, but does not by itself give $r_\xi=l$.
The rank of $\delta$ depends on the extension class, so the identity
\eqref{eq:tau-filtration-identity} need not be constant when only
$M$, $N$, and $l$ are fixed.

\noindent The term $\dim\Hom(M,N\otimes I_\xi)$ has a different role:
by Serre duality it equals $\dim\Ext^2(N\otimes I_\xi,M)$,
and hence contributes to the extension-stack dimension in
Remark~\ref{morfp}. It also bounds the tangent dimension of the
Quot fiber when the filtration is forgotten. It is not the rank of
$\delta$ and cannot replace that rank in
\eqref{eq:tau-filtration-identity}. On locally closed strata where
these dimensions and $\tau$ are constant and base change holds,
the extension contribution and the Grassmann contribution
$2(g-3-d+\tau)$ can be combined. Substituting
\eqref{eq:tau-filtration-bound} then gives the refined upper bound
in the parameter calculation below.

\begin{Remark} Consider the rational map
    $$h_{N,l}: \GG_{N,l}\rightarrow \WW^1_d(|L|),$$
    mapping a general point $(E,\Lambda)$ to the pair $(C_\Lambda, A_\Lambda)$, where $C_\Lambda$ is the degeneracy locus of the evaluation map $\text{ev}_\Lambda : \Lambda \otimes \OO_S \rightarrow E$, which is injective for a general $\Lambda \in G(2, H^0(S, E))$ and 
the fiber of $h_{N,l}$ at $(C,A)$ is the quotient stack
$$[\mathbb{P}(\Hom(E_{(C,A)},\omega_C\otimes A^\vee))^\circ/\text{Aut}(E_{(C,A)})].$$
As in \cite{lelli2013stability} any fibre of $h_{N,l}$ has dimension $-1$. 
\end{Remark}

\begin{Proposition}
    Assume that $\PP_{N,l}$ is non-empty and let $W$ be an irreducible component of $W_{N,l}$. Then
$$\dim (W) \leq  g+d+8-k,$$
where $k$ is the gonality of any curve in $|L|_s$.
\end{Proposition}
\begin{proof}
The slope ordering of $M$ and $N$ gives
$$\dim\Hom (M,N\otimes I_\xi)=\begin{cases}
    1,\ M\cong N, \xi=\emptyset\\
    0,\ \text{ otherwise}
\end{cases}$$
Indeed, strict inequality of slopes forces $\Hom(M,N\otimes I_\xi)=0$.
For equal slopes, a nonzero morphism $M\to N$ defines an effective
divisor of zero intersection with the ample class $L$, so $M\simeq N$.
In that case $\Hom(M,N\otimes I_\xi)=H^0(S,I_\xi)$, which has
dimension one for $l=0$ and zero for $l>0$.
It follows from Remark \ref{morfp} that the dimension of the fibers of $p$ is constant and equals $\langle v(N\otimes I_\xi),v(M)\rangle$,
unless $ M\cong N$ and $l = 0$, in which case it is $\langle v(N\otimes I_\xi),v(M)\rangle + 1$ then
$$\dim\EE_{N,l}=2l+\dim\Hom(M,N\otimes I_\xi)+\langle v(N\otimes I_\xi),v(M)\rangle.$$

\noindent It is well known (\cite[Proposition 2.2.8]{huybrechts2010geometry}) that the fiber of $q$ identifies with $\Quot_S(E,P)$ and that given an element $[\phi:E\rightarrow N\otimes I_\xi]$ then
$$\dim_{[\varphi]}\Quot_S(E,P)\leq \dim \Hom (K,N\otimes I_\xi)$$
where $K=\ker\varphi$ and the fiber is zero-dimensional when $\dim\Hom(M,N\otimes I_\xi)=0$; the exceptional case $M\simeq N$, $l=0$ must be treated separately. From Remark \ref{morfq} 
$\dim\PP_{N,l}\leq\dim\EE_{N,l}$ and if $\PP_{N,l}\neq \emptyset$
\begingroup
\begin{align*}
 \dim\GG_{N,l}
 &\leq2l+\langle v(N\otimes I_\xi),v(M)\rangle\\
 &\qquad+\dim\Hom(M,N\otimes I_\xi)+2(g-3-d+\tau)\\
 &=3l+2c_1(M)\cdot c_1(N)-(g-1)\\
 &\qquad+\dim\Hom(M,N\otimes I_\xi)+2(g-3-d+\tau)\\
 &=g+d-5-c_1(M)\cdot c_1(N)+2\tau
       +\dim\Hom(M,N\otimes I_\xi)\\
 &\leq g+d-5-c_1(M)\cdot c_1(N)\\
 &\qquad+2\min\{4,h^1(S,M)+h^1(S,N)+l-r_\xi\}\\
 &\qquad+\dim\Hom(M,N\otimes I_\xi)\\
 &\leq g+d+8-k.
\end{align*}
The final inequality assumes $\tau\leq4$, $\dim\Hom(M,N\otimes I_\xi)\leq1$, and
$c_1(M)\cdot c_1(N)\geq k-4$.
\endgroup
\noindent where we have used that $l=d-c_1(M) \cdot c_1(N)$, $c_1(E)=c_1(M)+c_1(N)$ and $c_1(E)^2=2g-2$. 
\end{proof}

\begin{Theorem}
Assume that general curves in $|L|_s$ have Clifford dimension $1$ and maximal gonality $k = 
\lfloor\frac{g+3}{2}\rfloor$.
 If $\rho(g, 1, d)>20$, any dominating component of $W^1_d(|L|)$ corresponds to $\mu_L$-stable LM bundles. In particular, if $C \in  |L|_s $ is general, the variety $W^1_d(|L|)$ is reduced and has the expected dimension $\rho(g, 1, d)$.

 \end{Theorem}
 \begin{proof}
\noindent From the above Proposition the dimension of any non-stable component $W$  of $\mathcal{W}_d^1(|L|)$ is bounded by
$$\dim W \le g + d +8 - k \leq g + d +8 - \frac{g+2}{2}=\frac{g}{2}+d+7.$$
On the other hand, any dominating component of $\mathcal{W}_d^1(|L|)$ must have dimension at least
$$\dim W \geq \dim |L|_s + \rho(g, 1, d) = g - 2 + \rho(g, 1, d) = g - 2 + (2d - g - 2) = 2d - 4.$$

\noindent Then
$$2d - 4 \leq \dim W \leq d + \frac{g}{2} +7,$$
and $d \leq \frac{g+2}{2} + 10$. However when $\rho>20$
$$d=\frac{\rho+g+2}{2}>\frac{g}{2} + 11.$$
 \end{proof}

 \section{Bridgeland Stability}

Let $S$ be a smooth complex projective abelian surface with
\begin{equation}\label{eq:picard-hypothesis}
 \NS(S)=\mathbb Z\cdot L,\qquad L\text{ ample}.
\end{equation}
For a smooth connected curve $i:C\hookrightarrow S$ in $|L|$, put
$g=1+L^2/2$. Let $A$ be globally generated of degree $d$, and write
\[
 U=H^0(C,A),\qquad n=r+1=\dim U,\qquad r>2,
 \qquad \chi_C=d+1-g.
\]
Thus $r$ denotes the dimension of the complete $g_d^r$, not the bundle
rank. Set $E=E_{C,A}$, $F=E^\vee$, $Q=i_*A$, and
$Q_B=i_*B$, where $B=\omega_C\otimes A^{-1}$. The defining sequences
are \eqref{eq:1} and \eqref{eq:2}. By
Theorem~\ref{thm:basic-lm-invariants},
\begin{equation}\label{eq:lm-vectors}
 v(Q)=(0,L,\chi_C),\quad v(Q_B)=(0,L,-\chi_C),\quad
 v(E)=(n,L,-\chi_C),\quad v(F[1])=(-n,L,\chi_C).
\end{equation}
Assume $\chi_C\neq0$ and choose the negative integer
$\chi=-|\chi_C|$. Following \cite[Definition 3.2]{bayer2018brill}, set
\begin{equation}\label{eq:classes}
 o=(1,0,0),\qquad
 R=\left\lfloor\frac{L^2}{-2\chi}\right\rfloor,\qquad
 w_k=(k-R,L,\chi).
\end{equation}
The plane $W_\chi=\operatorname{span}_{\mathbb R}\{o,w_0\}$ determines
the relevant wall. Let $\sigma_0$ be a geometric Bridgeland stability
condition on it and let $\WallCat$ be the abelian category of
$\sigma_0$-semistable objects of the common phase of $\OS$ and $w_k$.
Here geometric conditions are constructed on the tilted hearts
$\Coh^\beta S$; their existence and the required wall-and-chamber
properties are recalled in \cite[Theorem 2.1 and Proposition 2.2]{bayer2018brill},
with the large-volume identifications in \cite[Proposition 2.4]{bayer2018brill}.

Choose a nearby generic $\sigma_+$ with
$\nu_+(\OS)<\nu_+(w_k)$. For a fixed finite range containing all
classes used below, $\sigma_+$ can be chosen simultaneously in the
ordinary Gieseker chamber for nonnegative rank and the dual Gieseker
chamber for negative rank \cite[Proposition 3.1]{bayer2018brill}. In particular, stability of $F[1]=E^\vee[1]$
in the latter chamber is equivalent to Gieseker stability of $E$.
The description of stable factors and the adjacent-chamber criterion
\cite[Lemma 3.3, Proposition 3.4 and Corollary 3.5]{bayer2018brill} give, for
$T\in\WallCat$ of class $w_k$,
\[
 T\text{ is }\sigma_+\text{-stable}
 \quad\Longleftrightarrow\quad
 \Hom(T,P)=0\quad\text{for all }P\in\Pic^0(S).
\]
All exact sequences on the wall below are taken in $\WallCat$, not
necessarily in $\Coh(S)$.

For $T\in D^b(S)$, write
$H^j(T)=\Hom(\OS,T[j])$, $h^j(T)=\dim H^j(T)$, and set
\[
 M_k^h=\{T\in M_{\sigma_+}^s(w_k):h^0(T)=h\}.
\]
The number $h$ is a hypercohomology dimension and not necesarily equal to $n$.
For example, $h^0(F[1])=h^1(S,E)$ by Serre duality.
Bayer-Li's extension calculation gives
$\dim\Ext^1(T,\OS)=h-\chi$ for $T\in M_k^h$
\cite[Lemma 4.1]{bayer2018brill}.

\begin{Lemma}[Bayer-Li, Lemma 4.2]\label{lem:main}
Let $T\in M_k^h$.
\begin{enumerate}
\item[\textup{(a)}] For $W\subseteq H^0(T)$, evaluation gives an exact
sequence in $\WallCat$
\[
 0\longrightarrow\OS\otimes W\longrightarrow T
 \longrightarrow T'\longrightarrow0.
\]
Then $T'$ is $\sigma_+$-stable, has class $w_{k-\dim W}$, and
$h^0(T')\geq\dim W+\chi$.
\item[\textup{(b)}] For $V\subseteq\Ext^1(T,\OS)$, the extension defined
by the inclusion of $V$ gives an exact sequence
\[
 0\longrightarrow\OS\otimes V^\vee\longrightarrow\widetilde T
 \longrightarrow T\longrightarrow0.
\]
Then $\widetilde T$ is $\sigma_+$-stable, has class $w_{k+\dim V}$,
and $h^0(\widetilde T)\geq\dim V$.
\end{enumerate}
\end{Lemma}
For the proof, see \cite[Lemma 4.2]{bayer2018brill}.

The reduced object is defined by the unique exact sequence
\[
 0\longrightarrow J\longrightarrow T\longrightarrow T_{\red}
 \longrightarrow0,\qquad \Hom(\OS,T_{\red})=0,
\]
where $J$ is an iterated extension of copies of $\OS$
\cite[Lemma 3.6 and Definition 3.7]{bayer2018brill}. It is obtained by repeated
evaluation quotients; one quotient need not remove all sections since
$h^1(S,\OS)=2$. The construction excludes $\OS$ as a subobject of
$T_{\red}$, not as an arbitrary Jordan-H\"older factor.

\begin{Corollary}[Bayer-Li, Corollary 4.3]\label{cor:main}
For $T\in M_k^h$, its reduction $T_{\red}$ is $\sigma_+$-stable and
belongs to $M_{k_{\red}}^0$ for some $0\leq k_{\red}\leq k$.
\end{Corollary}
For the proof, see \cite[Corollary 4.3]{bayer2018brill}.

\begin{Theorem}[LM stability via Bayer-Li]
\label{higher-rank-lm-stability}
Let $S$ be a smooth complex projective abelian surface with
$\NS(S)=\mathbb Z\cdot L$, where $L$ is ample. Let $C\in|L|$ be smooth
and connected, put $g=1+L^2/2$, and let $A$ be a globally generated line
bundle of degree $d$ on $C$ with $h^0(C,A)=r+1$ and $r>2$.
If $d\neq g-1$, then $E_{C,A}$ is $\mu_L$-stable and
$L$-Gieseker-stable.
Moreover, for the adjacent condition $\sigma_+$ described above,
\[
 T=\begin{cases}
 E_{C,A}^\vee[1],&d<g-1,\\
 E_{C,A},&d>g-1
 \end{cases}
\]
is $\sigma_+$-stable, and $T_{\red}$ is $\sigma_+$-stable with
$H^0(T_{\red})=0$.
\end{Theorem}
\begin{proof}
Both $Q$ and $Q_B$ are Gieseker-stable: a proper nonzero subsheaf of
a line bundle on the integral curve $C$ has a zero-dimensional quotient
of positive length and hence smaller reduced Hilbert polynomial.
The chamber identifications above therefore apply to the sheaf with
negative Euler characteristic.

If $d<g-1$, then $Q\in M_R^n$. Rotating the evaluation sequence
\eqref{eq:1} gives the exact sequence in $\WallCat$
\[
 0\longrightarrow\OS\otimes U\longrightarrow Q
 \longrightarrow F[1]\longrightarrow0.
\]
Lemma~\ref{lem:main}(a) implies that $F[1]=E^\vee[1]$ is
$\sigma_+$-stable of class $w_{R-n}$. The dual Gieseker chamber then
gives Gieseker stability of $E$.

If $d>g-1$, then $Q_B$ has class $w_R$ and is $\sigma_+$-stable.
Serre duality and adjunction identify
\[
 \Ext^1(Q_B,\OS)\simeq H^1(C,B)^\vee\simeq H^0(C,A)=U.
\]
Under this identification, the class of the dual evaluation sequence
\eqref{eq:2} corresponds to the identity on $U$, as follows by dualizing
\eqref{eq:1}. Thus Lemma~\ref{lem:main}(b), applied to the full space
$V=U$, gives $\sigma_+$-stability of $E$, of class $w_{R+n}$.
The ordinary Gieseker chamber gives its Gieseker stability.

In either case $E$ is slope-semistable. If a proper saturated subsheaf
of rank $a$, $0<a<n$, had the same slope, its first Chern class would
be $mL$ for an integer $m$, and equality of slopes would give
$mL^2/a=L^2/n$, hence $mn=a$, which is impossible.
Therefore $E$ is $\mu_L$-stable. Finally,
Corollary~\ref{cor:main} applied to the stable object $T$ gives the
assertion about its reduction.
\end{proof}

The preceding theorem establishes stability of the Lazarsfeld-Mukai
bundle on $S$. The following corollary combines that result with
Feyzbakhsh's effective restriction criterion, see \cite{feyzbakhsh2022effective}. We first recall the normalized discriminant and then state Proposition~4.6 in the notation of the present setting.

Let $G$ be a $\mu_L$-stable reflexive sheaf of rank $q>1$ on $S$. We use Feyzbakhsh's normalized discriminant
\[
 \widetilde\Delta(G):=\left(\frac{c_1(G)\cdot L}{\operatorname{rk}(G)L^2}\right)^2-\frac{2\operatorname{ch}_2(G)}{\operatorname{rk}(G)L^2}.
\]
\begin{Theorem}[\cite{feyzbakhsh2022effective}, Prop.~4.6]
\label{feyz-prop}
Let $(S,L)$ be a smooth polarized projective surface such that $L^2$ divides $L\cdot L'$ for every $L'\in\Pic(S)$. Let $G$ be a $\mu_L$-stable reflexive sheaf of rank $q>1$ satisfying
\[
 \gcd\bigl(L^2\operatorname{rk}(G),c_1(G)\cdot L\bigr)=L^2.
\]
Define
\[
 m(G)_\pm=\min\left\{m\in\mathbb N:m'L^2\operatorname{rk}(G)-m\,c_1(G)\cdot L=\pm L^2\text{ for some }m'\in\mathbb Z\right\},
 \]
 \[
  \delta(G)=\frac{1}{q}\min\left\{\frac{1}{m(G)_+},\frac{1}{m(G)_-}\right\}.
\]
For any irreducible divisor $D\in|kL|$, the restriction $G|_D$ is $\mu$-stable if
\[
 k>\frac{\widetilde\Delta(G)}{\delta(G)}+\delta(G),
\]
and is $\mu$-semistable if the strict inequality is replaced by $\geq$. In particular, $G|_D$ is $\mu$-stable if
\[
 k>q(q-1)\widetilde\Delta(G)+\frac{1}{q(q-1)},
\]
and $\mu$-semistable if $>$ is replaced by $\geq$.
\end{Theorem}

In the application below, the Picard-rank-one hypothesis ensures the divisibility condition.

\begin{Corollary}
Use the notation of Theorem~\ref{higher-rank-lm-stability}, and assume $d<g-1$. Thus $U=H^0(C,A)$, $n=r+1$, $E=E_{C,A}$, and $F=E^\vee$. If $r>2$, then $F|_C$ is $\mu_{L|_C}$-stable. 
\end{Corollary}
\begin{proof}
The preceding theorem gives $\mu_L$-stability of $E$, and hence of its dual $F$, which is locally free. Recall that $\chi_C=d+1-g$ and $\operatorname{ch}(F)=(n,-L,-\chi_C)$. Using $L^2=2(g-1)$, the normalized discriminant becomes
\[
 \widetilde\Delta(F)=\left(\frac{c_1(F)\cdot L}{nL^2}\right)^2-\frac{2\operatorname{ch}_2(F)}{nL^2}=\frac{1}{r+1}\left(\frac{d}{g-1}-\frac{r}{r+1}\right).
\]
Because $\NS(S)=\mathbb Z\cdot L$, every divisor class is an integer multiple of $L$, so the divisibility condition holds. For $F$, the other hypothesis is satisfied since $\gcd(L^2\operatorname{rk}(F),c_1(F)\cdot L)=\gcd(nL^2,-L^2)=L^2$. Thus Proposition 4.6 applies to $C\in|L|$ with $k=1$, and its strict stability condition is
\[
 1>n(n-1)\widetilde\Delta(F)+\frac{1}{n(n-1)}.
\]
Substituting $n=r+1$ into the criterion gives
\[
 n(n-1)\widetilde\Delta(F)+\frac{1}{n(n-1)}=\frac{rd}{g-1}+\frac{1-r^3}{r(r+1)}<r+\frac{1-r^3}{r(r+1)}=\frac{r^2+1}{r(r+1)}<1,
\]
where the strict inequality uses $d<g-1$. Thus Theorem \ref{feyz-prop} implies that $F|_C$ is $\mu_{L|_C}$-stable.
\end{proof}

\bibliographystyle{plainurl}
\bibliography{Biblio}

\end{document}